\documentclass{jaums}
\usepackage{amsfonts,amsmath,amsthm,amssymb}

\def\gp#1{\langle #1 \rangle}
\def\m1{^{-1}}
\newcommand{\Y}{\mathbin{\mathsf Y}}

\theoremstyle{thmit} 
\newtheorem{thm}{Theorem}[section]
\newtheorem{cor}[thm]{Corollary}

\theoremstyle{thmrm} 

\newtheorem*{oldproof}{Proof}
\renewenvironment{proof}[1][{}]{\begin{oldproof}[#1]}{\qed\end{oldproof}}

\title[Powerful group of units]{Group algebras whose group of units is powerful}
\author{Victor Bovdi}
\address{Institute of Mathematics, University of Debrecen, \\
H-4010 Debrecen, P.O.B. 12, \\
Institute of Mathematics and Informatics,\\
College of Ny\'\i regyh\'aza, S\'ost\'oi \'ut 31/b, H-4410 Ny\'\i regyh\'aza,\\
Hungary}

\thanks{The research was supported by OTKA  No.K68383}



\begin{document}

\maketitle

\begin{abstract}
A $p$-group is called powerful if every commutator is a product of $p$\,th powers when $p$ is odd and a  product of fourth  powers when $p=2$. In the group algebra of a group $G$ of $p$-power order over a finite field of characteristic $p$, the group of normalized units is always a $p$-group. We prove that it is never powerful except, of course, when $G$ is abelian.
\end{abstract}

Throughout  this note $G$ is a finite $p$-group and
$F$ is a finite field of characteristic  $p$.  Let
$V(FG)=\{\sum_{g\in G}\alpha_gg\in FG\mid \sum_{g\in G}\alpha_g=1\}$ be the group of normalized units of the group algebra $FG$. Clearly $V(FG)$ is a finite $p$-group of order
$$
|V(FG)|=|F|^{|G|-1}.
$$
A $p$-group is called powerful if every commutator is a product of $p$\,th powers when $p$ is odd and a product of fourth  powers when $p=2$. The notion of powerful groups was
introduced in \cite{5}  and it  plays  an important
role in the study of finite $p$-groups (for example, see
\cite{2},  \cite{4} and  \cite{7}). Our main result is the following.

\begin{thm} \label{T:1}
The group of  normalized units $V(FG)$ of the group algebra $FG$ of a group $G$ of $p$-power order over a finite field $F$ of characteristic $p$,  is never powerful except, of course, when $G$ is abelian.
\end{thm}

In view of  the fact that a pro-$p$ group is powerful if and only if it is the limit of finite
powerful groups, this has an immediate consequence.

\begin{cor}
The group of normalized  units $V(F[[G]])$ of the completed group algebra $F[[G]]$ of a pro-$p$ group $G$ over a finite field $F$ of characteristic $p$,  is never powerful except, of course, when $G$ is abelian.
\end{cor}

We denote by $\zeta(G)$ the center of $G$. We say that $G=A\Y B$ is a central product of its subgroups $A$
and $B$ if $A$ and $B$ commute elementwise  and $G=\gp{A,B}$,
provided also that $A\cap B$ is the center of (at least) one of
$A$ and $B$. If $H$ is a subgroup of $G$, then  by ${\mathfrak I}(H)$ we denote the ideal of
$FG$ generated by the elements $h-1$ where  $h\in H$. Set
$(a,b)=a\m1b\m1ab$, where $a,b\in G$. Denote by $|g|$ the
order of $g\in G$.
Put $\Omega_k(G)=\gp{u\in G\mid u^{p^k}=1}$  and
$\widehat{H}=\sum_{g\in H}g\in FG$. If $H\trianglelefteq G$ is a
normal subgroup of $G$, then  $FG/{\mathfrak I}(H)\cong F[G/H]$  and
\begin{equation}\label{E:1}
V(FG)/(1+ {\mathfrak I}(H))\cong V(F[G/H]).
\end{equation}
We freely use the fact that every quotient of a powerful group is powerful (Lemma 2.2(i) in \cite{2}).
\begin{proof}
We prove the theorem by assuming that counterexamples exist, considering one of minimal order, and deducing a contradiction. Suppose then that $G$ is a counterexample of minimal order.  If $G$ had a nonabelian proper factor group $G/H$, that would be a smaller counterexample, for, by (\ref{E:1}),\; $V(F[G/H])$ would be a homomorphic image of the powerful group $V(FG)$. Thus all proper factor groups of $G$ are abelian, that is, $G$ is just nilpotent-of-class-$2$ in the sense of Newman \cite{6}. As Newman noted in the lead-up to his Theorem 1, this means
that the derived group has order $p$ and the centre is cyclic. Of course it follows that all $p$\,th powers are central, so the Frattini subgroup $\Phi(G)$ is central and also cyclic.

Suppose $p>2$. Then a finite $p$-group with only one subgroup of order $p$ is cyclic (Theorem 12.5.2 in \cite{3}), so $G$ must have a non-central subgroup $B=\langle b\rangle$ of order $p$. Now $(b,a)=c\neq1$ for some $a$ in $G$ and some $c$ in $G'$. Of course $\langle c\rangle=G'\leq\zeta(G)$,\;  $a^{-1}b^ia=b^ic^i=c^ib^i$\;  and\;  $b^i\widehat B=\widehat B$ for all $i$, so
\begin{equation}\label{E:2}
\begin{split}
(a\widehat B)^2
&=a^2(1+a^{-1}ba+\dots+a^{-1}b^{p-1}a)\widehat B\\
&=a^2(1+cb+\dots+c^{p-1}b^{p-1})\widehat B\\
&=a^2\widehat{G'}\widehat B.
\end{split}
\end{equation}
Noting that
\begin{equation}\label{E:3}
(\widehat{G'})^2=0,
\end{equation}
we get
\begin{equation}\label{E:4}
\begin{split}
(a\widehat B)^3=a^2\widehat{G'}\widehat B\cdot a\widehat B
&=a^2\widehat{G'}a^{-1}\cdot(a\widehat B)^2\\
&=a^2\widehat{G'}a^{-1}\cdot a^2\widehat{G'}\widehat B
=a^3(\widehat{G'})^2\widehat B=0.
\end{split}
\end{equation}
Therefore $|1+a\widehat B|=p$. We know from 4.12 of \cite{7} that $\Omega_1(V(FG))$ has exponent $p$, so we must have $((1+a\widehat B)b)^p=1$ as well. However,
\begin{equation}\label{E:5}
b^iab^{-i}=a(a,b^{-i})=ac^i=c^ia
\end{equation}
allows one to calculate that
\begin{align*}
((1+a\widehat B)b)^p
&=(1+a\widehat B)(1+bab^{-1}\widehat B)\cdots(1+b^{p-1}ab^{-(p-1)}\widehat B)\cdot b^p
&&\\
&=(1+a\widehat B)(1+ca\widehat B)\cdots(1+c^{p-1}a\widehat B)
&&\text{by (\ref{E:5})}\\
&=1+\widehat{G'}(a\widehat B)+\textstyle\frac12(p-1)\widehat{G'}(a\widehat B)^2
&&\text{by (\ref{E:4})}\\
&=1+\widehat{G'}(a\widehat B)+\textstyle\frac12(p-1)(\widehat{G'})^2a^2\widehat B
&&\text{by (\ref{E:2})}\\
&=1+ \widehat{G'}(a\widehat B)
&&\text{by (\ref{E:3})}\\
&\neq1.
\end{align*}
(To see that the third line equals the second, it helps to think in terms of polynomials with $a\widehat B$ as the indeterminate and $FG'$ as the coefficient ring, the critical point being that in the third line the coefficients of all positive powers of $a\widehat B$ are integer multiples of $\widehat{G'}$.) This contradiction completes the proof when $p>2$.

Next, we turn to the case $p=2$. Then $G'=\gp{c\mid c^2=1}$ and the ideal $\mathfrak I(G')$ is spanned by the elements of the form $\widehat{G'}g$, while $FG$ is spanned by the elements  $h$  of $G$. It is clear that $\widehat{G'}g$ and $h$ commute, because
$$
\widehat{G'}gh=\widehat{G'}(ghg^{-1}h^{-1})hg\qquad  \text{and}\qquad  \widehat{G'}(ghg^{-1}h^{-1})=\widehat{G'},
$$
so $\mathfrak I(G')$ is central in $FG$ and $1+\mathfrak I(G')$ is central in $V(FG)$. As $(\widehat{G'})^2=0$, it also follows $(\mathfrak I(G'))^2=0$ and so every element of $1+\mathfrak I(G')$ squares to $1$. As $V(FG)/(1+\mathfrak I(G'))\cong V(F[G/G'])$, the derived group $V'$ of $V(FG)$ lies in $1+\mathfrak I(G')$, a central subgroup of exponent 2. It follows that in $V(FG)$ all squares are central.

Let $w\in V'$. By Proposition 4.1.7 of \cite{5}, this is the fourth   power of some element $u$ of $V(FG)$. Write $u$ as $\sum_{g\in G}\alpha_gg$ with each $\alpha_g$ in $F$. In the commutative quotient modulo $\mathfrak I(G')$,\quad  $u^2=\sum_{g\in G}\alpha_g^2g^2$, hence
$$
u^2=v+\sum_{g\in G}\alpha_g^2g^2
$$
for some $v$ in $\mathfrak I(G')$. Of course then $v$ and all the $g^2$ are central in $FG$ and $v^2=0$, so we may conclude that\quad $w=u^4=\sum_{g\in G}\alpha_g^4g^4$. \;

In particular, as $V(FG)$ is not abelian, the exponent of $G$ must be larger than $4$. Recall that $\Phi(G)$ is central, the centre is cyclic, and $|G'|=2$, so Theorem 2 of \cite{1} applies and for this case gives the structure of $G$ as
$$
G=G_0\Y G_1\Y\cdots\Y G_r
$$
where $G_1,\dots,G_r$ are dihedral groups of order $8$ and $G_0$ is either cyclic of order at least $8$ (and in this case $r>0$) or an $M(2^{m+2})$ with $m>1$, where
$$
M(2^{m+2})=\langle\; a,b\; \mid\; a^{2^{m+1}}=b^2=1,\; a^b=a^{1+2^{m}}\;\rangle.
$$
One of the conclusions we need from this is that every fourth power in $G$ is already a fourth power in $G_0$, thus every element of $V'$ is an element of $FG_0^4$. In particular, when $w$ is the unique nontrivial element of $G'$, the linear independence of $G$ as subset of $FG$ implies that $w$ itself is the fourth power of some element of $G_0$.

It is easy to verify that, in $M(2^{m+2})$ with $m\geq1$, the inverse of the element $1+a+b$\;  is \; $(a^{2^m-3}+a^{-3}+a^{-2}+a^{-1})+(a^{2^m-2}+a^{2^m-2}+a^{-3})b$\;  and so
$$
(1+a+b,a)=(1+a^{2^m-2}+a^{-2})+(a^{2^m-2}+a^{2^m-1}+a^{-2}+a^{-1})b.
$$
Of course the left hand side is an element of $V'$, but the right hand side is not an element of $\langle a\rangle$. When $G_0\cong M(2^{m+2})$, this shows that there is an element in $V'$ which does not lie in $FG_0^4$. When $G_0$ is cyclic, then $G_1\cong M(2^{m+2})$ with $m=1$, and we have an element in $V'$ which does not even lie in $FG_0$. In either case, we have reached the promised contradiction and the proof of the Theorem is complete.
\end{proof}
The author would like to express his gratitude to L.~G.~Kov\'acs and particularly to the referee, for valuable remarks.

\bibliographystyle{abbrv}
\bibliography{Jams-paper}

\begin{thebibliography}{1}

\bibitem{1}
T.~R. Berger, L.~G. Kov{\'a}cs, and M.~F. Newman.
\newblock Groups of prime power order with cyclic {F}rattini subgroup.
\newblock {\em Nederl. Akad. Wetensch. Indag. Math.}, 42(1):13--18, 1980.

\bibitem{2}
J.~D. Dixon, M.~P.~F. du~Sautoy, A.~Mann, and D.~Segal.
\newblock {\em Analytic pro-{$p$} groups}, volume~61 of {\em Cambridge Studies
  in Advanced Math.}
\newblock Cambridge University Press, Cambridge, 1999.

\bibitem{3}
M.~Hall, Jr.
\newblock {\em The theory of groups}.
\newblock The Macmillan Co., New York, N.Y., 1959.

\bibitem{4}
L.~H{\'e}thelyi and L.~L{\'e}vai.
\newblock On elements of order {$p$} in powerful {$p$}-groups.
\newblock {\em J. Algebra}, 270(1):1--6, 2003.

\bibitem{5}
A.~Lubotzky and A.~Mann.
\newblock Powerful {$p$}-groups. {I}. {F}inite groups.
\newblock {\em J. Algebra}, 105(2):484--505, 1987.

\bibitem{6}
M.~F. Newman.
\newblock On a class of nilpotent groups.
\newblock {\em Proc. London Math. Soc. (3)}, 10:365--375, 1960.

\bibitem{7}
L.~Wilson.
\newblock On the power structure of powerful {$p$}-groups.
\newblock {\em J. Group Theory}, 5(2):129--144, 2002.

\end{thebibliography}

\setlength{\parindent}{0pt}

\end{document}